\documentclass[10pt]{article}
\usepackage{amsmath}
\usepackage{amsfonts}
\usepackage{mathrsfs}
\usepackage{amsfonts}
\usepackage{amssymb,amsfonts,amsmath, psfrag,eepic,colordvi,epsfig}
\topskip 
\numberwithin{equation}{section}
\newtheorem{theorem}{Theorem}[section]

\newtheorem{corollary}[theorem]{Corollary}

\newtheorem{conjecture}[theorem]{Conjecture}

\begin{document}
	\parskip 7pt
	
	\pagenumbering{arabic}
	\def\sof{\hfill\rule{2mm}{2mm}}
	\def\ls{\leq}
	\def\gs{\geq}
	\def\SS{\mathcal S}
	\def\qq{{\bold q}}
	\def\MM{\mathcal M}
	\def\TT{\mathcal T}
	\def\EE{\mathcal E}
	\def\lsp{\mbox{lsp}}
	\def\rsp{\mbox{rsp}}
	\def\pf{\noindent {\it Proof.} }
	\def\mp{\mbox{pyramid}}
	\def\mb{\mbox{block}}
	\def\mc{\mbox{cross}}
	\def\qed{\hfill \rule{4pt}{7pt}}
	\def\pf{\noindent {\it Proof.} }
	\textheight=21cm
	
	\begin{center}
		{\Large\bf   Parity results 
		for  the reciprocals  of 
		 false theta functions		}
	\end{center}
	
	\begin{center}

		Jing Jin$^{1}$,
		 Huan Xu$^{2}$ and Olivia X.M. Yao$^{3}$
		\\[10pt]
		
		$^{1}$College of Agricultural Information,\\
		Jiangsu Agri-animal
		Husbandry Vocational College,
		\\
		Taizhou, 225300,  Jiangsu,  P. R. China

		$^{2,3}$School of Mathematical Sciences, \\
		Suzhou University of Science and
		Technology, \\
		Suzhou,  215009, Jiangsu,
		P. R. China

		Email:
		$^{1}$jinjing19841@126.com,
		$^{2}$xu$\_$huanz@163.com, 
		$^{3}$yaoxiangmei@163.com

	\end{center}
	
	\noindent {\bf Abstract.}
	 Recently, Keith   investigated arithmetic properties  
	     for the 
	 reciprocals
	 of  some false theta functions and posed several conjectures. 
	 In this paper, we     prove some parity results 
	  for the 
	  reciprocals
	  of  some false theta functions by using some 
	   identities on Ramanujan's
	    general theta function. In particular, our results 
	     imply some conjectures of Keith.

		

	\section{Introduction}
	
	\allowdisplaybreaks
	
	The aim of this paper is to investigate 
	 parity results for the  
	 reciprocals
	 of  some false theta functions and confirm some 
	  conjectures given by Keith \cite{Keith} by utilizing 
	   several identities on  
	Ramanujan's
	general  theta
	function. Let us begin with some definitions. 
	 
	  Ramanujan's general theta function 
	 is  defined by
	\begin{align}\label{1-1}
		f(a,b)=\sum_{n=-\infty}^\infty
		a^{n(n+1)/2}b^{n(n-1)/2}, \qquad |ab|<1.
	\end{align}
	In Ramanujan's notation,
	the Jacobi triple product
	identity can be rewritten as 
	\begin{align}\label{1-2}
		f(a,b)=f(b,a)=(-a,-b,ab;ab)_\infty.
	\end{align}
	Here and throughout this paper,
	we use the following
	notation:
	\begin{align*}
		(a;q)_\infty &\;:=\prod_{k=0}^\infty
		(1-aq^k), \qquad
		(a;q)_n  :=\frac{(a;q)_\infty}
		{(aq^n;q)_\infty},\\
		(a_1,a_2,\ldots,a_k;q)_\infty
		&\;:=(a_1;q)_\infty (c_{2};q)_\infty
		\cdots (a_k;q)_\infty,
	\end{align*}
	where $q$ is a complex number
	with $|q|<1$. In addition, for all positive
	integers $m$, define
	\begin{align*}
		f_m:= (q^m;q^m)_\infty.
	\end{align*}
	
	Three important  special  cases
	of \eqref{1-1} are defined  by,  
	\begin{align}
		f( -q,-q^2)&\;=\sum_{n=-\infty}^\infty
		(-1)^n q^{n(3n-1)/2}=f_1,\label{1-3}\\
		f(-q,-q)&\;=\sum_{n=-\infty}^\infty
		(-1)^n q^{n^2}=\frac{f_1^2}{f_2},\label{1-4}\\
		f(-q,-q^3)&\;=\sum_{n=0}^\infty
		(-q)^{n(n+1)/2}
		=\frac{f_1f_4}{f_2}. \label{1-5}
	\end{align}
	Note that the reciprocals of the right
	sides
	of \eqref{1-3}--\eqref{1-5}
	are the generating functions of
	the ordinary partition function $p(n)$ \cite{Andrews-1976},
	the  overpartition
	function $\overline{p}
	(n)$
	\cite{Corteel},
	and the partition
	with odd parts distinct
	function \cite{Hirschhorn-1}, respectively,
	which are three of the most important
	types of partition functions. A number of
	congruences for the three 	partition functions have been proved;
	see for example
	\cite{HS-2005b,Hirschhorn-1,Ramanujan-1}.
	
	In 1917, Rogers \cite{Rogers} introduced false theta functions,
	which are series
	that would be classical theta functions
	except for changes in signs
	of an infinite number of terms. False theta functions closely
	resemble ordinary theta functions; however, they do not have the
	modular transformation properties that theta functions have.
	In his  notebooks \cite{Ramanujan-2}
	and lost notebook \cite{Ramanujan-3}, Ramanujan
	recorded many  
	identities involving  false theta functions
	$\Psi(a,b)$ which is defined by
	\begin{align}\label{1-6}
		\Psi(a,b):=\sum_{n=0}^\infty a^{n(n+1)/2} b^{n(n-1)/2}
		-\sum_{n=-\infty}^{-1}
		a^{n(n+1)/2}b^{n(n-1)/2}.
	\end{align}

	In recent years, there has been an increase in
	interest in false theta functions. They appear in a variety of
	contexts including modular and quantum modular forms,
	meromorphic Jacobi forms,
	quantum knot invariants,
	and many more; see for example \cite{Andrews-1981,Andrews-2007,Bringmann,Burson,Jennings-Shaffer,Wang}. 
	
	Very recently, Keith \cite{Keith}  studies 
	reciprocals
	of false theta functions
	and proved some interesting results
	on congruences,
	asymptotic bounds,
	and combinatorial identities. In particular, 
	  he showed 
	that for $n\geq 0$,
	\begin{align*}
		c_{5,1}(8n+5) &\; \equiv 0 \pmod 2,  \\
		c_{5,1}(32n+31) &\; \equiv 0 \pmod 4,
	\end{align*}
	where $c_{r,s}(n)$ is defined by
	\begin{align}\label{1-7}
		\sum_{n=0}^\infty c_{r,s}(n)q^n:=\frac{1}{\Psi(-q^r,q^s)}.
	\end{align}
	In addition, he  proved that $c_{9,1}(8n+4)\equiv 0 
	\pmod 2$
	 if $n$ cannot be represented in the form
	$n=10k^2-4k$ for integer $k$. 
	
	At the end of his paper, Keith \cite{Keith}
	posed the following conjecture
	on congruences modulo 2 for $c_{13,1}(n)$:
	
	\begin{conjecture} \cite{Keith}
		\label{conj-1}
		For $n\geq 0$,
		\begin{align}
			c_{13,1}(32n+23) &\;\equiv 0 \pmod 2,\label{1-8}\\[6pt]
			c_{13,1}(64n+63) &\;\equiv 0 \pmod 2,\label{1-9}\\[6pt]
			c_{13,1}(72n+j) &\;\equiv 0 \pmod 2,\label{1-10} 
		\end{align}
		where   $j\in\{15,21,39,69\}$.
	\end{conjecture}
	
	In this paper,  using some identities 
	 on Rmanujan's general theta function $f(a,b)$,
	  we prove some parity 
	 results 
	  for $c_{r,s}(n)$ where
	   $(r,s)\in \{(9,1),(7,3), (13,1),(11,3),(9,5)\}$. 
	  In particular,
	 we establish several infinite families 
	  of congruences modulo 2 which imply the congruences
	  given in Conjecture \ref{conj-1}.

	\section{Parity results for $c_{9,1}(n)$}	
	 
	 \begin{theorem}\label{Th-2-1}
	Let $p$ be  a prime with $p\equiv 
	11, \ 19 \pmod {20}$. For $n,k\geq 0$ with $p\nmid n$,
	\begin{align}\label{2-1}
	c_{9,1}\left(8 p^{2k+1}n+\frac{6p^{2k+2}+4}{5}\right)
	\equiv 0  \pmod 2 
	\end{align}
	and 
	\begin{align}\label{2-2}
	c_{9,1}\left(8p^{2k+1}n+\frac{26p^{2k+2}+4}{5}\right) \equiv 
	0 \pmod 2.
	\end{align}
\end{theorem}

\noindent{\it Proof.}
Using \eqref{1-2} and  the fact that
\[
(1-q)^2\equiv (1-q^2) \pmod 2,
\]
one can prove that 
 for all positive integers $r$ and $s$,
 \begin{align}\label{2-3}
 f(q^r,q^s)^2& \equiv f(q^{2r},q^{2s}) \pmod 2. 
 \end{align}
In view of \eqref{1-6} and \eqref{1-7}, we see  that 
  for  all positive integers $r$ and $s$,
 \begin{align}
 	\sum_{n=0}^\infty 
 	 c_{r,s}(n)q^n&=\frac{1}{\sum\limits_{n=0}^\infty (-1)^{n(n+1)/2}
 	  q^{(r+s)n^2/2+(r-s)n/2}-\sum\limits_{n=-\infty}^{-1} (-1)^{n(n+1)/2}
 	  q^{(r+s)n^2/2+(r-s)n/2}}
 	  \nonumber\\
 	  &\equiv  \frac{1}{\sum\limits_{n=-\infty}^\infty (-1)^{n(n+1)/2}
 	  	q^{(r+s)n^2/2+(r-s)n/2}}
 	  \nonumber\\
 	  &\equiv  \frac{1}{f(q^r,q^s)}
 	  \nonumber\\
 	  &\equiv  \frac{f(q^s,q^r)}{f(q^{2s},q^{2r})}
 	  \pmod 2.  \qquad \qquad ({\rm by}\ \eqref{2-3})
 	  \label{2-4}
 \end{align}
  It follows from \cite[Entry 30 (ii) and (iii)]{Berndt-1} that 
\begin{align}\label{2-5}
f(a,b)=f(a^3b,ab^3)+af(b/a,a^5b^3).
\end{align}
 Taking  $a=q$ and $b=q^9$ in \eqref{2-5} yields 
 \begin{align}\label{2-7}
 f(q,q^9)=f(q^{12},q^{28})+qf(q^8,q^{32}).
 \end{align}
 Setting $(r,s)=(9,1)$ in  \eqref{2-4}, we obtain 
 \begin{align}\label{2-8}
 	\sum_{n=0}^\infty 
 	c_{9,1}(n)q^n 
 	&\equiv  \frac{f(q,q^9)}{f(q^{2},q^{18})}
 	\pmod 2. 
 \end{align} 
Substituting \eqref{2-7}
into \eqref{2-8} and
extracting those
terms in which the power of  $q
$ is congruent to 0 modulo 2, then
  replacing
$q^{2}$ by $q$, we obtain
 \begin{align}\label{2-9}
 	\sum_{n=0}^\infty 
 	c_{9,1}(2n)q^n 
 	&\equiv  \frac{f(q^6,q^{14})}{f(q,q^{9})}
 	\nonumber\\
 	&\equiv  \frac{f(q,q^9)f(q^6,q^{14})}{f(q^2,q^{18})}
 	\pmod 2.  \qquad \qquad ({\rm by}\  \eqref{2-3})
 \end{align} 
 If we substitute  \eqref{2-7}
 into \eqref{2-9} and
  pick  out  those
 terms in which the power of  $q
 $ is congruent to 1 modulo 2, then
 divided  by  $q$  and replace 
 $q^{2}$ by $q$, we deduce that 
 \begin{align} \label{2-10}
 	\sum_{n=0}^\infty 
 	c_{9,1}(4n+2)q^n  
 	&\equiv  \frac{f(q^4,q^{16})f(q^3,q^{7})}{f(q,q^{9})}
 	\nonumber\\
 	&\equiv \frac{f(q^4,q^{16})}{f(q^2,q^{18})} f(q,q^9)f(q^3,q^7)
 	\pmod 2.    \qquad \qquad ({\rm by}\  \eqref{2-3})
 \end{align} 
 It follows from \cite[Entry 29]{Berndt-1} that 
 if $ab=cd$, then 
 \begin{align}\label{2-6}
 	f(a,b)f(c,d)=f(ac,bd)f(ad,bc)+af(b/c,ac^2d)f(b/d,acd^2).
 \end{align}
 Taking $(a,b,c,d)=(q,q^9,q^3,q^7)$
  in \eqref{2-6}, we arrive at 
  \begin{align}\label{2-11}
  f(q,q^9)f(q^3,q^7)=f(q^4,q^{16})f(q^8,q^{12})
  +qf(q^2,q^{18})f(q^6,q^{14}). 
  \end{align}
  If we  substitute  \eqref{2-11}
  into \eqref{2-10} and
  extract  those
  terms in which the power of  $q
  $ is congruent to 0 modulo 2, then
   replace 
  $q^{2}$ by $q$, we arrive at 
  \begin{align}\label{2-12}
  	\sum_{n=0}^\infty 
  	c_{9,1}(8n+2)q^n  
  	&\equiv \frac{f(q^2,q^{8})^2 f(q^4,q^6)}{f(q,q^{9})} 
  	\nonumber\\
  	 & \equiv  \frac{f(q,q^{4})^4 f(q^2,q^3)^2 }{f(q,q^{9})}  
  	\pmod 2.  \qquad \qquad ({\rm by}\  \eqref{2-3})
  \end{align} 
  It is easy to check
   that 
   \begin{align}\label{2-13}
       \frac{f(q,q^{4})^4 f(q^2,q^3)^2 }{f(q,q^{9})} 
   &=  \frac{f(q,q^{4})^3
   	  f(q^2,q^3)^2 f(q^4,q^6) f_5}{f_{10}^2} 
   	  \nonumber\\
   	  &\equiv
   	   \frac{f(q,q^{4})^3
   	   	f(q^2,q^3)^4  f_5}{ f_{10}^2} 
   \pmod 2.\qquad \qquad ({\rm by} \ \eqref{2-3})
   \end{align}
   Combining \eqref{2-12} and \eqref{2-13}
    yields 
    \begin{align}\label{2-14}
    	\sum_{n=0}^\infty 
    	c_{9,1}(8n+2)q^n  
    	&\equiv \frac{f(q,q^{4})^3
    		f(q^2,q^3)^4  f_5}{ f_{10}^2}   
    	\pmod 2.   
    \end{align}  
   It is easy to verify that 
   \begin{align}\label{2-15}
   f(q,q^4)f(q^2,q^3)=\frac{f_2f_5^3}{f_1f_{10}}
   .
   \end{align}
   In light of \eqref{2-14}
    and \eqref{2-15},
     \begin{align}
    	\sum_{n=0}^\infty
    	c_{9,1}(8n+2)q^n& \equiv \frac{f_2^3f_{5}^{10}}{f_1^3f_{10}^5}f(q^2,q^3) 
    	\nonumber\\
    	& \equiv \frac{f_2^2}{f_1}f(q^2,q^3) 
    	 \pmod 2, \label{2-16}
    \end{align}
    where we have used  the fact that 
    for any positive integer $k$,
   \begin{align}
   	f_k^2\equiv f_{2k} \pmod 2. \label{2-17}
   \end{align}
   It is well known that
   \begin{align}
   	\sum_{m=0}^\infty q^{m(m+1)/2}=\frac{f_2^2}{f_1}.
   	\label{2-18}
   \end{align}
   In view of \eqref{1-1}, \eqref{2-16} and \eqref{2-18}, 
   \begin{align*}
   	\sum_{n=0}^\infty
   	c_{9,1}(8n+2)q^n \equiv  \sum_{m=0}^\infty
   	 q^{m(m+1)/2} \sum_{k=-\infty}^\infty
   	  q^{(5k^2+k)/2} \pmod 2,
   \end{align*}
which can be rewritten as 
\begin{align*}
	\sum_{n=0}^\infty
	c_{9,1}(8n+2)q^{40n+6} \equiv  \sum_{m=0}^\infty
	\sum_{k=-\infty}^\infty
	q^{ 5(2m+1)^2 +(10k+1)^2 } \pmod 2.
\end{align*}
Therefore,  if $40n+6$
is not of the form
$ 5x^2+y^2$, then $c_{9,1}(8n+2)\equiv 0 \pmod 2$.
Note that if $N$ is
of the form $  5x^2+y^2$
and  $p$ is a prime
with $p\equiv 11,19 \pmod {20} $,
then $\nu_p(N)$ is even since
$\left(
\frac{-5  }{p}\right)=-1$. Here $\nu_p(N)$
denotes the highest power of $p$
dividing
$N$ and $\left(
\frac{ \cdot }{p}\right)$ denotes
the Legendre symbol.
It is easy to check that if $p\nmid n$,
then for any
$\alpha\in \mathbb{N}$,
\[
\nu_p\left(40\left(p^{2\alpha+
	1}n+\frac{3(p^{2\alpha +2}-1)
}{20}\right)+6\right)=\nu_p(40p^{2\alpha+1}n
+6p^{2\alpha+2
})=2\alpha+1
\]
is odd and  $40\left(p^{2\alpha+1}n+ \frac{3(p^{2\alpha+2 }-1)
}{20}\right)+6$ is not of the form
$ 5x^2+y^2$. Therefore,
if  $p$ is a prime
with $p\equiv 11,19 \pmod {20} $
and $p\nmid n$, then for any
$\alpha\in \mathbb{N}$,
\begin{align*}
	c_{9,1}\left(8\left(p^{2\alpha+
		1}n+\frac{3(p^{2\alpha+2 }-1) }{20}\right)+2\right) 
	\equiv 0 \pmod 2,
\end{align*}
from which, we arrive at  \eqref{2-1}.

  Substituting \eqref{2-11}
  into \eqref{2-10} and
  picking out  those
  terms in which the power of  $q
  $ is congruent to 1 modulo 2, then
  dividing   by  $q$  and replacing
  $q^{2}$ by $q$, we get 
  \begin{align}
  	\sum_{n=0}^\infty 
  	c_{9,1}(8n+6)q^n  
  	&\equiv f(q^2,q^{8})f(q^3,q^7) 
  	\pmod 2.  \label{2-19}
  \end{align}  
  In view of \eqref{1-1} and \eqref{2-19}
\begin{align}
\sum_{n=0}^\infty 
c_{9,1}(8n+6)q^n 
& \equiv \sum_{m,k=-\infty}^\infty
 q^{5m^2+3m+5k^2+2k} \pmod 2. \label{2-20}
\end{align}
It follows from \eqref{2-20} that 
\begin{align*}
	\sum_{n=0}^\infty 
	c_{9,1}(8n+6)q^{20n+13} &  \equiv \sum_{m,k=-\infty}^\infty
	q^{(10m+3)^2 +(10k+2)^2} \pmod 2.
\end{align*}
Thus,  if $20n+13$
is not of the form
$ x^2+y^2$, then $c_{9,1}(8n+6)$ is even.
If  $N$ is
of the form $  x^2+y^2$
and  $p$ is a prime
with $p\equiv 11,19 \pmod {20} $,
then $\nu_p(N)$ is even because 
$\left(
\frac{-1  }{p}\right)=-1$.  It is easy to check that if $p\nmid n$,
then for any
$\alpha\in \mathbb{N}$,
\[
\nu_p\left(20\left(p^{2\alpha+
	1}n+\frac{13(p^{2\alpha +2}-1)
}{20}\right)+13\right)=\nu_p(20p^{2\alpha+1}n
+13p^{2\alpha+2
})=2\alpha+1
\]
is odd and  $20\left(p^{2\alpha+1}n+ \frac{13(p^{2\alpha+2 }-1)
}{24}\right)+13$ is not of the form
$ x^2+y^2$. Thus,
if  $p$ is a prime
with $p\equiv 11,19 \pmod {20} $
and $p\nmid n$, then for any
$\alpha\in \mathbb{N}$,
\begin{align*}
	c_{9,1}\left(8\left(p^{2\alpha+
		1}n+\frac{13(p^{2\alpha+2 }-1) }{20}\right)+6\right) 
		\equiv 0 \pmod 2,
\end{align*}
which yields \eqref{2-2}. This completes
 the proof of Theorem \ref{Th-2-1}.  \qed

	\section{Parity results for $c_{7,3}(n)$}

	 	\begin{theorem}\label{Th-3-1}
	 		(1) 	For $n\geq 0$,
	 		\begin{align}\label{3-1}
	 		c_{7,3}(16n+5)\equiv 0\pmod 2.
	 		\end{align}
	 		(2) For $n\geq 0$, 	$c_{7,3}(16n+13)
	 		\equiv 0 \pmod 2$
	 		if and only if 
	 		$n$ can not be represented in the form 
	 		$n=5k^2+4k$ for some integer $k$. 
	 		\\
	 		(3) 	Let $p$ be  a prime with $p\equiv 
	 		11, \ 19 \pmod {20}$. For $n,k\geq 0$ with $p\nmid n$,
	 		\begin{align}\label{3-2}
	 		c_{7,3}\left(8p^{2k+1}n+\frac{34p^{2k+2}+1}{5}\right) \equiv 
	 		0 \pmod 2.
	 		\end{align}
	 	\end{theorem}
	 	
	 	{\bf Remark}.  Note that the second part (2) of Theorem 
	 	\ref{Th-3-1} is analogous to the parity results 
	 	for $c_{9,1}(16n+14)$ proved by Keith \cite{Keith}. 
	 	
	 	\noindent{\it Proof}.
	 	Setting $(r,s)=(7,3)$ in \eqref{2-4} yields 
	 	\begin{align}
	 	\sum_{n=0}^\infty c_{7,3}(n)q^n \equiv \frac{f(q^3,q^7)}{f(q^6,q^{14})}
	 	\pmod 2. \label{3-3}
	 	\end{align}
	 	Taking $(a,b)=(q^3,q^7)$ in \eqref{2-5} yields 
	 		\begin{align}\label{3-4}
	 	f(q^3,q^7)=f(q^{16},q^{24})+q^3f(q^4,q^{36}).
	 	\end{align}
	 	Substituting \eqref{3-4}
	 	into \eqref{3-3} and
	 	extracting those
	 	terms in which the power of  $q
	 	$ is congruent to 1 modulo 2, then
	 	dividing   by  $q$  and replacing
	 	$q^{2}$ by $q$, we obtain
	 	\begin{align}\label{3-5}
	 	\sum_{n=0}^\infty c_{7,3}(2n+1)q^n& \equiv  q \frac{f(q^2,q^{18})}{f(q^3,q^{7})} \nonumber\\
	 	&\equiv  q \frac{f(q^2,q^{18})}{f(q^6,q^{14})}f(q^3,q^7)
	 	\pmod 2. \qquad ({\rm by}\ \eqref{2-3})
	 	\end{align}
	 Substituting \eqref{3-4} into \eqref{3-5},
	  we obtain 
	 \begin{align*}
	 	\sum_{n=0}^\infty c_{7,3}(2n+1)q^n \equiv  q \frac{f(q^2,q^{18})}{f(q^6,q^{14})}(f(q^{16},q^{24})+q^3f(q^4,q^{36}))
	 	\pmod 2.  
	 \end{align*}
	 Thus,we can immediately read off the following:
	 		\begin{align}\label{3-6}
	 	\sum_{n=0}^\infty c_{7,3}(4n+1)q^n &\equiv  q^2 \frac{f(q,q^{9})}{f(q^3,q^{7})}f(q^2,q^{18})
	 	\nonumber\\
	 	& \equiv  q^2 \frac{f(q^2,q^{18})
	 	}{f(q^6,q^{14})}  f(q,q^{9})f(q^3,q^{7})
	 	\pmod 2 \qquad \qquad ({\rm by}\ \eqref{2-3})
	 	\end{align}
	 	and 
	 		\begin{align}\label{3-7}
	 	\sum_{n=0}^\infty c_{7,3}(4n+3)q^n &\equiv   \frac{f(q^8,q^{12})}{f(q^3,q^{7})}  f(q,q^{9})
	 	\nonumber\\
	 	& \equiv   \frac{f(q^8,q^{12})}{f(q^6,q^{14})}  f(q,q^{9})f(q^3,q^{7}) 
	 	\pmod 2. \qquad \qquad ({\rm by}\ \eqref{2-3})
	 	\end{align}
	 	If we substitute \eqref{2-11}
	 	into \eqref{3-6} and
	 	pick out  those
	 	terms in which the power of  $q
	 	$ is congruent to 1 modulo 2, then
	 	divided  by  $q$  and replace 
	 	$q^{2}$ by $q$, we obtain
	 	\begin{align*}
	 	\sum_{n=0}^\infty c_{7,3}(8n+5)q^n & \equiv    qf(q,q^{9})^2 \nonumber\\
	 	& \equiv qf(q^2,q^{18})
	 	\pmod 2, \qquad ({\rm by}\ \eqref{2-3})
	 	\end{align*}
	 	which yields \eqref{3-1}
	 	 and 
	 	\[
	 		 	\sum_{n=0}^\infty c_{7,3}(16n+13)q^n \equiv    f(q,q^{9})
	 		 	\pmod 2. 
	 	\]
	 	 In addition,
	 	 it follows from \eqref{1-1} and the above congruence that 
	 	\[
	 	\sum_{n=0}^\infty c_{7,3}(16n+13)q^n \equiv   \sum_{k=-\infty}^\infty q^{5k^2+4k}
	 	\pmod 2.
	 	\]
The second part (2) of Theorem \ref{Th-3-1}
 follows from the above congruence. 
	 	
	 	Substituting \eqref{2-11}
	 	into \eqref{3-7} and
	 	picking out  those
	 	terms in which the power of  $q
	 	$ is congruent to 1 modulo 2, then
	 	dividing   y  $q$  and replacing
	 	$q^{2}$ by $q$, we obtain
	 	\begin{align}\label{3-8}
	 	\sum_{n=0}^\infty
	 	c_{7,3}(8n+7) q^n 
	 	\equiv f(q,q^9)f(q^4,q^6) \pmod 2.
	 	\end{align}
	 	In light of \eqref{1-1}
	 	 and \eqref{3-8}, 
	 	\begin{align}\label{3-9}
	 	\sum_{n=0}^\infty
	 	c_{7,3}(8n+7) q^n 
	 	\equiv  \sum_{m,k=-\infty}^\infty q^{5k^2+4k+5m^2+m} \pmod 2.
	 	\end{align}
	 	We can rewrite \eqref{3-9} as 
	 	\begin{align}\label{3-10}
	 	\sum_{n=0}^\infty
	 	c_{7,3}(8n+7) q^{20n+17} 
	 	\equiv  \sum_{m,k=-\infty}^\infty q^{(10k+4)^2+(10m+1)^2} \pmod 2.
	 	\end{align}
Using \eqref{3-10} and the same method for proving  \eqref{2-2}, we can get \eqref{3-2}.
 This completes the proof of Theorem \ref{Th-3-1}. \qed

	 	\section{Parity results for $c_{13,1}(n)$}
	 	
	 	In this section, we prove some parity results
	 	 for $c_{13,1}(n)$. In particular, we confirm Conjecture 
	 	  \ref{conj-1}. 
	 	
	 	 The following theorem generalizes \eqref{1-8} and \eqref{1-9}. 
	 	\begin{theorem}
	 		\label{Th-4-1}
	 		For $n,k\geq 0$, 
	 		\begin{align}
	 			c_{13,1}\left(2^{3k+5}n+\frac{19\times 2^{3k+3}+9}{7}\right) \equiv &\; 0 
	 			\pmod 2, \label{4-1}\\
	 			c_{13,1}\left(2^{3k+6}n+\frac{27\times 2^{3k+4}+9}{7}\right) \equiv &\;  0 
	 			\pmod 2, \label{4-2}\\
	 			c_{13,1}\left(2^{3k+7}n+\frac{3\times 2^{3k+5}+9}{7}\right) \equiv &\;  0 
	 			\pmod 2. \label{4-3}
	 		\end{align}
	 		
	 	\end{theorem}
	 	
	 	Setting $k=0$ in \eqref{4-1} and \eqref{4-2},
	 	 we obtain the following corollary.
	 	
	 	\begin{corollary}
	 		Congrunces \eqref{1-8} and \eqref{1-9} are true. 
	 	\end{corollary}

	 	\noindent{\it Proof.}
	 	 Taking $(r,s)=(13,1)$ in \eqref{2-4} yields 
	 	 \begin{align}\label{4-4}
	 	 \sum_{n=0}^\infty
	 	 c_{13,1}(n)q^n  \equiv \frac{f(q,q^{13})}{f(q^2,q^{26})} \pmod 2.
	 	 \end{align}
	 	 Setting $(a,b)=(q,q^{13})$ in \eqref{2-5}, we arrive at 
	 	 	 	 \begin{align}\label{4-5}
	 	 f(q,q^{13})=f(q^{16},q^{40})+qf(q^{12},q^{44}).
\end{align}
Substituting \eqref{4-5}
into \eqref{4-4} and
extracting those
terms in which the power of  $q
$ is congruent to 1 modulo 2, then
dividing    by  $q$  and replacing
$q^{2}$ by $q$, we deduce that 
	 	 \begin{align}\label{4-6}
	 	 \sum_{n=0}^\infty c_{13,1}(2n+1)q^n&
	 	 \equiv\frac{f(q^{6},q^{22})}{
	 	 f(q,q^{13})}
	 	 \nonumber\\
	 	 & \equiv \frac{f(q^6,q^{22})}{f(q^2,q^{26})}f(q,q^{13})
	 	  \pmod 2.  \qquad ({\rm by}\  \eqref{2-3})
	 	 \end{align}
	 	 Substituting \eqref{4-5}
	 	 into \eqref{4-6} and
	 	 picking out  those
	 	 terms in which the power of  $q
	 	 $ is congruent to 1 modulo 2, then
	 	 dividing   by  $q$  and replacing
	 	 $q^{2}$ by $q$, we obtain 
	 	 \begin{align}\label{4-7}
	 	 \sum_{n=0}^\infty c_{13,1}(4n+3)q^n &
	 	  \equiv \frac{f(q^3,q^{11})}{f(q,q^{13})}f(q^6,q^{22})
	 	  \nonumber\\
	 	  &
	 	 \equiv \frac{ f(q^6,q^{22}) }{f(q^2,q^{26})} 
	 	 f(q,q^{13})f(q^3,q^{11}) 
	 	 \pmod 2.  \qquad ({\rm by}\  \eqref{2-3})
	 	 \end{align}
	 	 Setting $(a,b,c,d)=(q,q^{13},q^3,q^{11})$
	 	  in \eqref{2-6} yields  
	 	 \begin{align}
	 	 f(q,q^{13})f(q^3,q^{11})
	 	 =f(q^4,q^{24})f(q^{12},q^{16})+qf(q^{10},q^{18})f(q^2,q^{26}).
	 	 \label{4-8}
	 	 \end{align}
 If we substitute   \eqref{4-8}
into \eqref{4-7} and
pick out  those
terms in which the power of  $q
$ is congruent to 1 modulo 2, then
divided   by  $q$  and replace 
$q^{2}$ by $q$, we have 
\begin{align}\label{4-9}
	 	 \sum_{n=0}^\infty
	 	 c_{13,1}(8n+7)q^n\equiv f(q^3,q^{11})f(q^5,q^9)
	 	 \pmod 2.
	 	 \end{align}
	 	 Putting $(a,b,c,d)=(q^3,q^{11},q^5,q^9)$
 in \eqref{2-6} yields 
  \begin{align}\label{4-10}
	 	 f(q^3,q^{11})f(q^5,q^9)=f(q^8,q^{20})f(q^{12},q^{16})
	 	 +q^3f(q^6,q^{22})f(q^2,q^{26}).
\end{align}
In light of  \eqref{4-9} and \eqref{4-10},
 \begin{align}\label{4-11}
 	\sum_{n=0}^\infty
 	c_{13,1}(8n+7)q^n\equiv f(q^8,q^{20})f(q^{12},q^{16})
 	+q^3f(q^6,q^{22})f(q^2,q^{26})
 	\pmod 2,
 \end{align}
 which   yields 
\begin{align}\label{4-12}
	 	 \sum_{n=0}^\infty
	 	 c_{13,1}(16n+7)q^n \equiv f(q^4,q^{10})f(q^6,q^8)
	 	  \pmod 2,
\end{align}
	 	 and 
\begin{align}\label{4-13}
	 	 \sum_{n=0}^\infty
	 	 c_{13,1}(16n+15)q^n \equiv q  f(q,q^{13}) f(q^3,q^{11})
	 	  \pmod 2.
	 	 \end{align}
It follows from \eqref{4-12} that for $n\geq 0$, 
	 	 \begin{align}\label{4-14}
	 	 c_{13,1}(32n+23)\equiv 0 \pmod 2. 
	 	 \end{align}
Thanks to 	\eqref{4-8} and 
\eqref{4-13},
 \begin{align}\label{4-15}
 	\sum_{n=0}^\infty
 	c_{13,1}(16n+15)q^n \equiv q f(q^4,q^{24})f(q^{12},q^{16})+q^2f(q^{10},q^{18})f(q^2,q^{26})
 	\pmod 2,
 \end{align}
	 which yields 	
	 	\begin{align}\label{4-16}
	 	\sum_{n=0}^\infty
	 	c_{13,1}(32n+15)q^n \equiv q  f(q,q^{13}) f(q^5,q^{9})
	 	\pmod 2
	 	\end{align}
	 	and 
	 	\begin{align}\label{4-17}
	 	\sum_{n=0}^\infty
	 	c_{13,1}(32n+31)q^n \equiv f(q^2,q^{12})f(q^{6},q^{8})
	 	\pmod 2.
	 	\end{align}
It follows from \eqref{4-17} that
 for $n\geq 0$, 
\begin{align}\label{4-18}
c_{13,1}(64n+63)\equiv 0 \pmod 2.
\end{align}
Setting $(a,b,c,d)=(q,q^{13},q^5,q^9)$
 in \eqref{2-6} yields 
\begin{align}\label{4-19}
f(q,q^{13})f(q^5,q^9)=f(q^6,q^{22})f(q^{10},q^{18})+qf(q^8,q^{20})f(q^4,q^{24}). 
\end{align}
In view of \eqref{4-16} and \eqref{4-19}, 
	\begin{align}\label{4-20}
	\sum_{n=0}^\infty
	c_{13,1}(32n+15)q^n \equiv q  f(q^6,q^{22})f(q^{10},q^{18})+q^2f(q^8,q^{20})f(q^4,q^{24})
	\pmod 2, 
\end{align}
from which, we arrive at 
\begin{align}\label{4-21}
\sum_{n=0}^\infty
c_{13,1}(64n+15)q^n \equiv q  f(q^4,q^{10}) f(q^2,q^{12})
\pmod 2
\end{align}
and 
\begin{align}\label{4-22}
\sum_{n=0}^\infty
c_{13,1}(64n+47)q^n \equiv  f(q^3,q^{11}) f(q^5,q^{9}) \pmod 2.
\end{align}
It follows from \eqref{4-21} that for $n\geq 0$,
\begin{align}\label{4-23}
c_{13,1}(128n+15)\equiv 0 \pmod 2. 
\end{align}
In light of \eqref{4-9}  and \eqref{4-22}, 
 \begin{align}\label{4-24}
c_{13,1}(64n+47) \equiv c_{13,1}(8n+7)
\pmod 2.
\end{align}
By \eqref{4-24}  and mathematical induction, 
\begin{align}\label{4-25}
c_{13,1}\left(2^{3k+3}n+\frac{5\times 2^{3k+3}+9}{7}\right) \equiv c_{13,1}(8n+7)
\pmod 2.
\end{align}
Replacing $n$ by $4n+2$, $8n+7$ and $16n+1$
 in \eqref{4-25}
  and using \eqref{4-14},
   \eqref{4-18} and \eqref{4-23}, we arrive at \eqref{4-1},
    \eqref{4-2} and \eqref{4-3}, respectively. 
     This completes the proof of Theorem \ref{Th-4-1}. \qed 

\begin{theorem}\label{Th-4-2}
	Let $p$ be   a prime with $p\equiv 
	15, 27 \pmod {28}$.  For $m,n,k\geq 0$ with $p\nmid n$,
	\begin{align}\label{4-26}
	c_{13,1}\left(2^{3k+3}p^{2m+1}n+\frac{5\times 
		 2^{3k+3}p^{2m+2}+9}{7}\right) 
	\equiv 0  \pmod 2.
	\end{align}
\end{theorem}

\noindent{\it Proof.} In view of \eqref{1-1} and \eqref{4-9}, 
\begin{align}\label{4-27}
\sum_{n=0}^\infty
c_{13,1}(8n+7)q^n\equiv \sum_{m,k=-\infty}^\infty
 q^{7m^2+4m+7k^2+2k} \pmod 2. 
\end{align}
Using \eqref{4-27} and the same method for proving \eqref{2-2}, we can prove that if 
  $p$ is  a prime with $p\equiv 
 15, 27 \pmod {28}$, then  for $m,n\geq 0$ with $p\nmid n$, 
 \begin{align}\label{4-28}
 c_{13,1}\left(8 p^{2m+1}n+\frac{40p^{2m+2}+9}{7}\right)
 \equiv 0  \pmod 2.
 \end{align}
 Replacing $n$ by $p^{2m+1}n+\frac{5(p^{2m+2}-1)}{7}\ (p\nmid n)$
  in \eqref{4-25} and utilizing \eqref{4-28}, we arrive at 
  \eqref{4-26}.  The proof is complete. \qed 
  
  	\begin{theorem}\label{Th-4-3}
  	Let $p$ be  a prime with $p\equiv 
  	15, \ 27 \pmod {28}$. For $n,k\geq 0$ with $p\nmid n$,
  	\begin{align}\label{4-29}
  	c_{13,1}\left(8p^{2k+1}n+\frac{26p^{2k+2}+9}{7}\right) \equiv 
  	0 \pmod 2.
  	\end{align}
  \end{theorem}
  
  \noindent{\it Proof.}
  Substituting \eqref{4-5}
  into \eqref{4-6} and
  picking out  those
  terms in which the power of  $q
  $ is congruent to 0 modulo 2, then
    replacing
  $q^{2}$ by $q$, we obtain  
  \begin{align}
  	\sum_{n=0}^\infty c_{13,1}(4n+1)q^n &
  	\equiv \frac{f(q^3,q^{11})}{f(q,q^{13})}f(q^8,q^{20})
  	\nonumber\\
  	&
  	\equiv \frac{ f(q^8,q^{20}) }{f(q^2,q^{26})} 
  	f(q,q^{13})f(q^3,q^{11}) 
  	\pmod 2.   \qquad ({\rm by}\  \eqref{2-3}) \label{4-30}
  \end{align}
  If we substitute   \eqref{4-8}
  into \eqref{4-30} and
  pick out  those
  terms in which the power of  $q
  $ is congruent to 1 modulo 2, then
  divided   by  $q$  and replace 
  $q^{2}$ by $q$, we deduce that 
\begin{align}
\sum_{n=0}^\infty
c_{13,1}(8n+5)q^n \equiv  f(q^4,q^{10}) f(q^5,q^{9})
\pmod 2.  \qquad ({\rm by} \ \eqref{2-3}) \label{4-31}
\end{align}
Based on \eqref{1-1} and \eqref{4-31},
\begin{align}
	\sum_{n=0}^\infty
	c_{13,1}(8n+5)q^n \equiv \sum_{m,k=-\infty}^\infty
	 q^{7m^2+3m+7k^2+2k} 	\pmod 2. \label{4-32}
\end{align}
Using \eqref{4-32} and the same method for proving \eqref{2-2}, we can prove Theorem \ref{Th-4-3}. This completes
 the proof. \qed  
 
 To conclude this section, we give a proof 
  of \eqref{1-10}.
  
  \begin{corollary}
  	Congruence \eqref{1-10} is true.
  \end{corollary}
 
 \noindent{\it Proof.}
   It is easy to check that  
\begin{align}\label{4-33}
	\sum_{m=-\infty}^\infty
	q^{7m^2+4m}=\sum_{m=-\infty}^\infty
	 q^{63m^2+12m}+\sum_{m=-\infty}^\infty
	 q^{63m^2+54m+11}+\sum_{m=-\infty}^\infty
	 q^{63m^2+96m+36}
\end{align}	
and 
\begin{align}\label{4-34}
\sum_{k=-\infty}^\infty
q^{7k^2+2k}=
\sum_{k=-\infty}^\infty
q^{63k^2+6k}+\sum_{k=-\infty}^\infty
q^{63k^2+48k+9}+\sum_{k=-\infty}^\infty
q^{63k^2+90k+32}.
\end{align}
  Substituting \eqref{4-33} and \eqref{4-34}
 into \eqref{4-27} and
 picking out  those
 terms in which the power of  $q
 $ is congruent to 1 modulo 3, then
  dividing by $q$ and 
 replacing
 $q^{3}$ by $q$, we obtain  
 \begin{align*}
 	\sum_{n=0}^\infty
 	c_{13,1}(24n+15)q^n\equiv 
 	 q^{14}\sum_{m,k=-\infty}^\infty q^{21m^2+18m+21k^2+30k}
 	\pmod 2,
 \end{align*}
 from which, we deduce that   that 
 for $n\geq 0$,
 \begin{align}\label{4-35}
 c_{13,1}(72n+15) \equiv c_{13,1}(72n+39) \equiv 0 \pmod 2. 
 \end{align}
In addition, 
\begin{align}\label{4-36}
	\sum_{m=-\infty}^\infty
	q^{7m^2+3m}=\sum_{m=-\infty}^\infty
	q^{63m^2+9m}+\sum_{m=-\infty}^\infty
	q^{63m^2+51m+10}+\sum_{m=-\infty}^\infty
	q^{63m^2+93m+34}. 
\end{align}	
 If we substitute   \eqref{4-34} and 
 \eqref{4-36}
 into \eqref{4-32} and
pick out  those
terms in which the power of  $q
$ is congruent to 2 modulo 3, then
divided   by  $q^2$  and replace 
$q^{3}$ by $q$, we obtain  
\begin{align*}
	\sum_{n=0}^\infty
	c_{13,1}(24n+21)q^n\equiv 
	q^{10}\sum_{m,k=-\infty}^\infty q^{21m^2+3m+21k^2+30k}
	\pmod 2,
\end{align*}
which implies that for $n\geq 0$,
\begin{align}\label{4-37}
c_{13,1}(72n+21)
\equiv c_{13,1}(72n+69) \equiv 0 \pmod 2.
\end{align}
Congruence \eqref{1-10} follows from \eqref{4-35} and 
\eqref{4-37}. 
 The proof is complete.  \qed 
	 	 
	 		\section{Parity results
	 		for $c_{11,3}(n)$}

	 	 	\begin{theorem}
	 		\label{Th-5-1}
	 		For $n,k\geq 0$, 
	 		\begin{align}
	 			c_{11,3}\left(2^{3k+5}n+\frac{3\times 2^{3k+3}+4}{7}\right) \equiv &\; 0 
	 			\pmod 2, \label{5-1}\\
	 			c_{11,3}\left(2^{3k+6}n+\frac{19\times 2^{3k+4}+4}{7}\right)\equiv &\;  0 
	 			\pmod 2, \label{5-2}\\
	 			c_{11,3}\left(2^{3k+7}n+\frac{27\times 2^{3k+5}+4}{7}\right) \equiv &\;  0 
	 			\pmod 2. \label{5-3}
	 		\end{align}
	 		
	 	\end{theorem}

	 	\noindent{\it Proof.} 
	 	 Setting  $(r,s)=(11,3)$ in \eqref{2-4}, we obtain 
	 	\begin{align}
	 		\sum_{n=0}^\infty
	 		c_{11,3}(n)q^n  \equiv \frac{f(q^3,q^{11})}{f(q^6,q^{22})} \pmod 2.  \label{5-4}
	 	\end{align}
	 	Taking  $(a,b)=(q^3,q^{11})$ in \eqref{2-5}
	 	 yields 
	 	\begin{align}
	 		f(q^3,q^{11})=f(q^{20},q^{36})+q^3 
	 		f(q^{8},q^{48}). \label{5-5}
	 	\end{align}
	 	Substituting \eqref{5-5}
	 	into \eqref{5-4} and
	 	extracting those
	 	terms in which the power of  $q
	 	$ is congruent to 0 modulo 2, then
	 	  replacing
	 	$q^{2}$ by $q$, we arrive  at 
	 	\begin{align}
	 		\sum_{n=0}^\infty c_{11,3}(2n)q^n&
	 		\equiv\frac{f(q^{10},q^{18})}{
	 			f(q^3,q^{11})}
	 		\nonumber\\
	 		& \equiv \frac{f(q^{10},q^{18})}{f(q^6,q^{22})}f(q^3,q^{11})
	 		\pmod 2.  \qquad ({\rm by}\  \eqref{2-3}) \label{5-6}
	 	\end{align}
	 		Substituting \eqref{5-5}
	 	into \eqref{5-6} and
	 	extracting those
	 	terms in which the power of  $q
	 	$ is congruent to 0 modulo 2, then
	 	replacing
	 	$q^{2}$ by $q$, we arrive at
	 	\begin{align}
	 		\sum_{n=0}^\infty c_{11,3}(4n)q^n 
	 		& \equiv \frac{f(q^{5},q^{9})}{f(q^3,q^{11})}f(q^{10},q^{18})
	 		\nonumber\\
	 		& \equiv \frac{f(q^{10},q^{18})}{f(q^6,q^{22})}
	 		 f(q^3,q^{11})f(q^5,q^9) 
	 		 \pmod 2.  \qquad ({\rm by} \ \eqref{2-3}) \label{5-7}
	 	\end{align}
	 	Substituting \eqref{4-10}
	 	into \eqref{5-7} and
	 	extracting those
	 	terms in which the power of  $q
	 	$ is congruent to 1 modulo 2, then
	 	dividing  them by  $q$  and replacing
	 	$q^{2}$ by $q$, we have 
	 	\begin{align}
	 		\sum_{n=0}^\infty c_{11,3}(8n+4)q^n 
	 	 	& \equiv qf(q,q^{13})f(q^5,q^9) 
	 		\pmod 2. \label{5-8}
	 	\end{align}
	 	It follows from \eqref{4-16}  and \eqref{5-8} that
	 	 for $n\geq 0$,
\begin{align}\label{5-9}
	 	 c_{11,3}(8n+4)\equiv c_{13,1}(32n+15) \pmod 2.
	 	 \end{align}
Replacing $n$ by $2^{3k+2}n+\frac{3(2^{3k}-1)}{7}$
 in \eqref{5-9} and utilizing \eqref{4-3}, we obtain \eqref{5-1}.	 
 
 Replacing $n$ by $2^{3k+3}n+\frac{19\times 
 	 2^{3k+1}-3}{7}$
 in \eqref{5-9}  yields 	 
	 \begin{align}
	 	c_{11,3}\left(2^{3k+6}n+\frac{19\times 2^{3k+4}+4}{7}\right)\equiv
	 		c_{13,1}\left(2^{3k+8}n+\frac{19\times 2^{3k+6}+9}{7}\right) \pmod 2. \label{5-10}
	 \end{align}	 
	 Replacing $k$ by $k+1$
	  in \eqref{4-1} yields  
	 \[
	 	c_{13,1}\left(2^{3k+8}n+\frac{19\times 2^{3k+6}+9}{7}\right) \equiv 0 \pmod 2,
	 \]
	 from which with \eqref{5-10}, we arrive at \eqref{5-2}. 
	 
	 Replacing $n$ by $2^{3k+4}n+\frac{27\times 2^{3k+2}-3}{7}$
	  in \eqref{5-9} yields 
	  \begin{align}
	  	c_{11,3}\left(2^{3k+7}n+\frac{27\times 2^{3k+5}+4}{7}\right) \equiv
	  	c_{13,1}\left(2^{3k+9}n+\frac{27\times 2^{3k+7}+9}{7}\right)  \pmod 2. \label{5-11}
	  \end{align}	
	  It follows from \eqref{4-2} that 
	  \[
	  c_{13,1}\left(2^{3k+9}n+\frac{27\times 2^{3k+7}+9}{7}\right)  
	  \equiv 0 \pmod 2.
	  \]
	  Congruence \eqref{5-3} follows from \eqref{5-11} and 
	   the above congruence.  
	 	  The proof of Theorem \ref{Th-5-1} is complete. \qed

	 	  \begin{theorem}\label{Th-5-2}
	 	  	Let $p$ be   a prime with $p\equiv 
	 	  	15, 27 \pmod {28}$.  For $m,n,k\geq 0$ with $p\nmid n$,
	 	  	\begin{align}
	 	  		c_{11,3}\left(2^{3k+3}p^{2m+1}n+\frac{3 \times 
	 	  			2^{3k+3}p^{2m+2}+4}{7}\right) 
	 	  		\equiv 0  \pmod 2. \label{5-12}
	 	  	\end{align}
	 	  \end{theorem}
	 	  
\noindent{\it Proof.}	  Replacing $n$ by $4n+1$ in \eqref{4-25} yields
\begin{align}
c_{13,1}(32n+15)\equiv c_{13,1}\left(2^{3k+5}n+\frac{3\times 2^{3k+5}+9}{7} \right) \pmod 2.  \label{5-13}
\end{align}
Replacing $n$ by $2^{3k}n+\frac{3(2^{3k}-1)}{7}$ in \eqref{5-9}, we get 
\begin{align} 
	c_{11,3}\left(2^{3k+3}n+\frac{3\times 2^{3k+3}+4}{7}\right) \equiv c_{13,1}\left(2^{3k+5}n+\frac{3\times 2^{3k+5}+9}{7}\right) 
	\pmod 2. \label{5-14}
\end{align}	 
Combining \eqref{5-9}, \eqref{5-13} and 
\eqref{5-14} yields 
\begin{align} 
	c_{11,3}\left(2^{3k+3}n+\frac{3\times 2^{3k+3}+4}{7}\right) \equiv c_{11,3}(8n+4) 
	\pmod 2. \label{5-15}
\end{align}	 
In view of \eqref{1-1} and \eqref{5-8}, 
	 	  \begin{align}
	 	  	\sum_{n=0}^\infty c_{11,3}(8n+4)q^{n-1} 
	 	  	& \equiv \sum_{m,k=-\infty}^\infty
	 	  	 q^{7m^2+6m+7k^2+2k}  
	 	  	\pmod 2. \label{5-16}
	 	  \end{align}
	 	  Using \eqref{5-16} and the same method for proving \eqref{2-2}, we can prove that if 
	 	  $p$ is  a prime with $p\equiv 
	 	  15, 27 \pmod {28}$, then  for $m,n\geq 0$ with $p\nmid n$, 
	 	  \begin{align}
	 	  	c_{11,3}\left(8 p^{2m+1}n+\frac{24p^{2m+2}+4}{7}\right)
	 	  	\equiv 0  \pmod 2. \label{5-17}
	 	  \end{align}
	 Replacing $n$ by $ p^{2m+1}n+\frac{3(p^{2m+2}-1)}{7}$	
	 in \eqref{5-15} and using \eqref{5-17}, we arrive at \eqref{5-12}.
	  This completes the proof of Theorem \ref{Th-5-2}.
	  \qed

	 	\begin{theorem}\label{Th-5-3}
	 		Let $p$ be  a prime with $p\equiv 
	 		15, \ 27 \pmod {28}$. For $n,k\geq 0$ with $p\nmid n$,
	 		\begin{align}\label{5-18}
	 		c_{11,3}\left(8p^{2k+1}n+\frac{122p^{2k+2}+4}{7}\right)
	 		 \equiv 
	 		0 \pmod 2.
	 		\end{align}
	 	\end{theorem}
	 	
	\noindent{\it Proof.}
	 	Substituting \eqref{5-5}
	 into  \eqref{5-6} and
	 extracting those
	 terms in which the power of  $q
	 $ is congruent to 1 modulo 2, then
	 dividing  by  $q$  and replacing
	 $q^{2}$ by $q$, we arrive at
		\begin{align}
		\sum_{n=0}^\infty c_{11,3}(4n+2)q^n 
		& \equiv q \frac{f(q^{5},q^{9})}{f(q^3,q^{11})}f(q^{4},q^{24})
		\nonumber\\
		& \equiv q \frac{f(q^{4},q^{24})}{f(q^6,q^{22})}
		f(q^3,q^{11})f(q^5,q^9) 
		\pmod 2.  \qquad ({\rm by} \ \eqref{2-3}) \label{5-19}
	\end{align}
	If we substitute \eqref{4-10}
	into \eqref{5-19} and
	 pick out  those
	terms in which the power of  $q
	$ is congruent to 0 modulo 2, then
	   replace 
	$q^{2}$ by $q$, we arrive at
	 	\begin{align}\label{5-20}
\sum_{n=0}^\infty
 	c_{11,3}(8n+2) q^n \equiv q^2 f(q^2,q^{12})f(q,q^{13})
 	 \pmod 2. 
	\end{align}
In light of \eqref{1-1} and \eqref{5-20},
	\[
\sum_{n=0}^\infty
c_{11,3}(8n+2) q^n \equiv \sum_{m,k=-\infty}^\infty 
 q^{7m^2+6m+7k^2+5k+2}  \pmod 2.
\]
 Utilizing  the above congruence 
  and the same method for proving \eqref{2-2}, we can 
  prove \eqref{5-18}. The proof of Theorem 
   \ref{Th-5-3} is complete. \qed 
	 	
	 	\section{Parity results
	 	  for $c_{9,5}(n)$}

	 		\begin{theorem}
	 		\label{Th-6-1}
	 		For $n,k\geq 0$, 
	 		\begin{align}
	 			c_{9,5}\left(2^{3k+5}n+\frac{27\times 2^{3k+3}+1}{7}\right) \equiv &\;  0 
	 			\pmod 2,  \label{6-1}\\
	 			c_{9,5}\left(2^{3k+6}n+\frac{3\times 2^{3k+4}+1}{7}\right) \equiv &\;  0 
	 			\pmod 2,\label{6-2}\\
	 				c_{9,5}\left(2^{3k+7}n+\frac{19\times 2^{3k+5}+1}{7}\right) \equiv &\; 0 
	 			\pmod 2. \label{6-3}
	 		\end{align}
	 		
	 	\end{theorem}

	\noindent{\it Proof.}
	 Setting $(r,s)=(9,5)$ in \eqref{2-4}
	 yields 
	 \begin{align}\label{6-4}
	 	\sum_{n=0}^\infty c_{9,5}(n)q^n
	 	 \equiv \frac{f(q^5,q^9)}{f(q^{10},q^{18})} \pmod 2.
	 \end{align}
	Taking $(a,b)=(q^5,q^9)$
	 in \eqref{2-5} yields 
	 \begin{align}\label{6-5}
	 f(q^5,q^9)=f(q^{24},q^{32})+q^5f(q^4,q^{52}).
	 \end{align}
	Substituting  \eqref{6-5}
 into  \eqref{6-4} and
extracting those
terms in which the power of  $q
$ is congruent to 1 modulo 2, then
dividing    by  $q$  and replacing
$q^{2}$ by $q$, we arrive at
\begin{align}\label{6-6}
	\sum_{n=0}^\infty c_{9,5}(2n+1)q^n
	&\equiv q^2 \frac{f(q^2,q^{26})}{f(q^{5},q^{9})} \nonumber\\
	&\equiv q^2 \frac{f(q^2,q^{26})}{f(q^{10},q^{18})} f(q^5,q^9)  \pmod 2.  \qquad ({\rm by}\  \eqref{2-3})
\end{align}
	If we substitute  \eqref{6-5} 
into  \eqref{6-6} and
extract  those
terms in which the power of  $q
$ is congruent to 1 modulo 2, then
divided    by  $q$  and replace 
$q^{2}$ by $q$, we arrive at
\begin{align}\label{6-7}
	\sum_{n=0}^\infty c_{9,5}(4n+3)q^n
	 &\equiv q^3 \frac{f(q,q^{13})}{f(q^{5},q^{9})} f(q^2,q^{26}) 
	 \nonumber\\
	 &\equiv q^3\frac{f(q^2,q^{26})}{
	 f(q^{10},q^{18})} f(q,q^{13})f(q^5,q^9) \pmod 2.  \qquad ({\rm by} \ \eqref{2-3})
\end{align}
Substituting \eqref{4-19} 
into  \eqref{6-7} and
extracting those
terms in which the power of  $q
$ is congruent to 1 modulo 2, then
dividing   by  $q$  and replacing
$q^{2}$ by $q$, we arrive at
\begin{align}\label{6-8}
	\sum_{n=0}^\infty c_{9,5}(8n+7)q^n
  &\equiv q f(q,q^{13})f(q^3,q^{11})  \pmod 2.
\end{align}
In light of \eqref{4-13} and \eqref{6-8},
\begin{align}\label{6-9}
c_{9,5}(8n+7)\equiv c_{13,1}(16n+15)  \pmod 2.
\end{align}
Replacing $n$ by $2^{3k+2}n+\frac{27\times 2^{3k}-6}{7}$
 in \eqref{6-9} and using \eqref{4-2}, we get \eqref{6-1}. 

Replacing $n$ by $2^{3k+3}n+\frac{3\times 2^{3k+1}-6}{7}$
in \eqref{6-9} and using \eqref{4-3}, we get \eqref{6-2}. 

Replacing $n$ by $2^{3k+4}n+\frac{19\times 
 2^{3k+2}-6}{7}$ in \eqref{6-9} yields 
 \begin{align}\label{6-10}
 c_{9,5}\left(2^{3k+7}n+\frac{19\times
  2^{3k+5}+1}{7}\right)\equiv c_{13,1}\left(2^{3k+8}n+\frac{19\times
   2^{3k+6}+9}{7}\right)  \pmod 2.
 \end{align}
It follows from \eqref{4-1} that 	
\begin{align}\label{6-11}
c_{13,1}\left(2^{3k+8}n+\frac{19\times
	2^{3k+6}+9}{7}\right) \equiv 0\pmod 2.
\end{align}
Congruence \eqref{6-3} follows from \eqref{6-10} and \eqref{6-11}. 
\qed 

\begin{theorem}\label{Th-6-2}
	Let $p$ be   a prime with $p\equiv 
	15, 27 \pmod {28}$.  For $m,n,k\geq 0$ with $p\nmid n$,
	\begin{align} \label{6-12}
		c_{9,5}\left(2^{3k+3}p^{2m+1}
		n+\frac{3\times 2^{3k+4}p^{2m+2}+1}{7}\right) 
		\equiv 0  \pmod 2.
	\end{align}
\end{theorem}

\noindent{\it Proof.}
 Replacing $n$ by $2n+1$ in \eqref{4-25}, we have
 \begin{align}\label{6-13}
 	c_{13,1}\left(2^{3k+4}n+\frac{3\times 2^{3k+5}+9}{7}\right)
 	\equiv c_{13,1}(16n+15) \pmod 2.
 \end{align}
  Replacing $n$ by $2^{3k}n+\frac{3\times 2^{3k+1}-6}{7}$ in 
  \eqref{6-9} yields
\begin{align} \label{6-14}
	c_{9,5}\left(2^{3k+3}n+\frac{3\times 2^{3k+4}+1}{7}\right) \equiv c_{13,1}\left(2^{3k+4}n+\frac{3\times 2^{3k+5}+9}{7}\right) 
	\pmod 2.
\end{align}	 
Combining \eqref{6-9}, \eqref{6-13} and \eqref{6-14}, we get 
\begin{align} 
	c_{9,5}\left(2^{3k+3}n+\frac{3\times 2^{3k+4}+1}{7}\right)  \equiv c_{9,5}(8n+7) 
	\pmod 2. \label{6-15}
\end{align}	 
In light of \eqref{1-1} and \eqref{6-8},
\begin{align}
	\sum_{n=0}^\infty c_{9,5}(8n+7)q^n
	&\equiv \sum_{m,k=-\infty}^\infty
	 q^{7m^2+6m+7k^2+4k+1}   \pmod 2. \label{6-16}
\end{align}
Using \eqref{6-16} and the same method 
 for proving \eqref{2-2},
  we can prove that if 
  $p$ is  a prime with $p\equiv 
  15, 27 \pmod {28}$, then  for $m,n\geq 0$ with $p\nmid n$, 
  \begin{align}
  	c_{9,5}\left(8 p^{2m+1}n+\frac{48p^{2m+2}+1}{7}\right)
  	\equiv 0  \pmod 2. \label{6-17}
  \end{align}
Replacing $n$ by $p^{2m+1}n+\frac{6(p^{2m+2}-1)}{7}$
 in \eqref{6-15} and using \eqref{6-17}, we get 
 \eqref{6-12}.
  This completes the proof
   of Theorem \ref{Th-6-2}. \qed  
	 	
	 	\begin{theorem} \label{Th-6-3}
	 		Let $p$ be  a prime with $p\equiv 
	 		15, \ 27 \pmod {28}$. For $n,k\geq 0$ with $p\nmid n$,
	 		\begin{align}\label{6-18}
	 		c_{9,5}\left(8p^{2k+1}n+\frac{34p^{2k+2}+1}{7}\right) \equiv 
	 		0 \pmod 2.
	 		\end{align}
	 	\end{theorem}

\noindent{\it Proof.}
	If we substitute  \eqref{6-5} 
into  \eqref{6-6} and
 pick out  those
terms in which the power of  $q
$ is congruent to 0 modulo 2, then
  replace 
$q^{2}$ by $q$, we arrive at
\begin{align} \label{6-19}
\sum_{n=0}^\infty
c_{9,5}(4n+1)q^n &
 \equiv q\frac{f(q,q^{13})}{f(q^5,q^9)}f(q^{12},q^{16})
 \nonumber\\
 &\equiv q\frac{f(q^{12},q^{16})}{f(q^{10},q^{18})}f(q,q^{13})f(q^5,q^9)
\pmod 2.  \qquad ({\rm by}\  \eqref{2-3})
\end{align}
Substituting \eqref{4-19}  
into  \eqref{6-19} and
extracting those
terms in which the power of  $q
$ is congruent to 1 modulo 2, then
  replacing
$q^{2}$ by $q$, we arrive at
	 	\[
	 	\sum_{n=0}^\infty
	 	c_{9,5}(8n+5)q^n \equiv f(q^6,q^8)f(q^3,q^{11})
	 	\pmod 2.
	 	\]
	Thanks to \eqref{1-1} and the above congruence,
\begin{align}\label{6-20}
\sum_{n=0}^\infty
c_{9,5}(8n+5)q^n \equiv \sum_{m,k=-\infty}^\infty
 q^{7m^2+m+7k^2+4k} 
\pmod 2.
\end{align}
Using \eqref{6-20} and the same method for proving \eqref{2-2},
 we can arrive at \eqref{6-18}. This completes
  the proof. 	   \qed 
 
 \section*{Statements and Declarations}
 
 \noindent{\bf Funding.}
 This work was supported by
 the National Natural Science Foundation of
 China  (grant no.
 12371334) and  the Qinlan project of
 Jiangsu province of China and the research project of Jiangsu Agri-animal Husbandry Vocational College (No. NSF2025CB22).
 
 \noindent{\bf Author Contributions.}
 The authors contributed equally to the preparation of this article. All authors read and
 approved the final manuscript.
 
 \noindent{\bf Competing Interests.}
 The authors declare that they have
 no conflict of interest.

 \noindent{\bf Data Availability Statements.} Data sharing not applicable to this
 article as no datasets were generated or analysed during the current
 study.

\end{document}